\documentclass[11pt]{article}%
\usepackage{amssymb}
\usepackage{latexsym}
\usepackage{amsmath}

\newtheorem{theorem}{Theorem}[section]

\newtheorem{corollary}[theorem]{Corollary}

\newtheorem{definition}[theorem]{Definition}

\newtheorem{lemma}[theorem]{Lemma}

\newtheorem{proposition}[theorem]{Proposition}
\newtheorem{remark}[theorem]{Remark}

\begin{document}

\begin{center}\Large A Modified Version of Free Orbit-Dimension\\ of von
Neumann Algebras
\end{center}

\centerline{Don Hadwin\hspace{4cm} Weihua Li}
 \centerline{Mathematics Department, University of New Hampshire, Durham, NH 03824}
 \centerline{don@unh.edu\hspace{4cm} whli@unh.edu}

\vspace{0.3cm} \noindent\textbf {Abstract}  Based on the notion of
free orbit-dimension introduced by D. Hadwin and J. Shen
\cite{Don-Shen}, we introduce a new invariant
 on finite von Neumann algebras that do not necessarily act on
separable Hilbert space. We show that this invariant is independent
on the generating set, and we extend some results in \cite{Don-Shen}
to von Neumann algebras that are not finitely generated.

\section{Introduction}
\hspace{1.5em}The theory of free entropy and free entropy dimension
was developed by D. Voiculescu in the 1990's and it is one of the
most powerful and exciting new tools in the theory of von Neumann
algebras. D. Voiculescu \cite{Vo1} \cite{Vo} introduced the concept
of free entropy in relation to his free probability theory and the
concept of free entropy dimension, and he used them to prove that
the free group factors do not contain Cartan subalgebras, which
answered a long-standing open problem. Later this was generalized by
L. Ge \cite{Ge}, who showed that the free group factors do not
contain a simple masa. L. Ge \cite{Ge1} used free entropy to give
the first example of a separable prime II$_1$ factor. Later, L. Ge
and J. Shen \cite{Ge-Shen} computed the free entropy dimension of
some II$_1$ factors with property T, including ${\cal L}(SL({\Bbb
Z},2m+1))\ (m\leq 1)$. Recently, D. Hadwin and J. Shen
\cite{Don-Shen} introduced a new invariant, the upper free
orbit-dimension of a finite von Neumann algebra, which is closely
related to Voiculescu's free entropy dimension. Using their new
invariant, they generalized and simplified the proofs of most of the
applications of free entropy dimension to finite von Neumann
algebras.

Here we introduce a new invariant, $\mathfrak{K}_3$, which is a
modification of the upper free orbit-dimension, $\mathfrak{K}_2$;
when $\mathfrak{K}_2$ is defined, $\mathfrak{K}_3=\infty\cdot
\mathfrak{K}_2$. We then extend the domain of $\mathfrak{K}_3$ to
all finite von Neumann algebras that can be embedded into some
ultrapower of the hyperfinite II$_1$ factors. This includes algebras
acting on nonseparable Hilbert spaces.

 The organization of the paper is as follows. In section 2, we recall the definition of free orbit-dimension,
 and introduce a new invariant $\mathfrak{K}_3$
on von Neumann algebras. In section 3, we prove:

(1)$\mathfrak{K}_3({\mathcal S})=\mathfrak{K}_3({\mathcal G})$ when
$W^*({\cal S})=W^*({\cal G})$ (independence of the generators),

 (2) if ${\mathcal N}_1\cap{\mathcal N}_2$
is diffuse, then $\mathfrak{K}_3\left(({\mathcal N}_1\cup {\mathcal
N}_2)''\right)\leq \mathfrak{K}_3({\mathcal
N}_1)+\mathfrak{K}_3({\mathcal N}_2),$

(3) $\mathfrak{K}_3\left(W^*({\cal N}\cup
\{y\})\right)\leq\mathfrak{K}_3({\mathcal N})$ whenever there exist
normal operators $a$ and $b$ in $\mathcal N$ without common
eigenvalues such that $ay=yb\neq 0$.

 In section 4 we prove:

(1) if $\{{\mathcal{M}}_{\lambda }\}_{\lambda \in
\Lambda }$ is a family of von Neumann algebras such that each $\mathcal{M}%
_{\lambda }$ has a central net of Haar untiaries, and $\alpha $ is a
nontrivial ultrafilter on $\Lambda $, then $
\mathfrak{K}_{3}(\prod^{\alpha }\mathcal{M}_{\lambda })=0, $

(2) if $\mathcal{M}$ is a von Neumann algebra with a central net of
Haar unitaries, then $\mathfrak{K}_{3}({\mathcal{M}})=0$,

(3) if $\mathbb{F}$ is a free group with the standard generating set
$G$ satisfying $|G|\geq 2$, then $\mathfrak{K}_{3}({\mathcal{L}}_{\mathbb{F}%
})=\infty$,

(4) if ${\mathcal N}_1$ and ${\mathcal N}_2$ are mutually commuting
diffuse subalgebras of $\mathcal M$, then $\mathfrak{K}_3(W^*({\cal
N}_1\cup {\cal N}_2))=0$,

(5) if $\cal M$ is a II$_1$ factor and $\mathfrak{K}_3({\cal
M})=\infty$, then $\cal M$ is prime (i.e., cannot be written as a
tensor product of two II$_1$ factors).

In section 5, we show how our invariant leads naturally to a
canonical decomposition of torsion-free groups into a union of
certain self-normalizing subgroups so that the intersection of any
two of them is $\{e\}$. We completely describe this decomposition
for free groups, and we present a related question for the free
group factor ${\cal L}_{{\Bbb F}_2}$.

All of the free entropy concepts require the von Neumann algebra
$\cal M$ under consideration can be tracially embedded into an
ultrapower of the hyperfinite II$_1$ factor. Throughout this paper,
we assume that all the von Neumann algebras we consider can be
embedded.
\section{Preliminaries}

\hspace{1.5em} First we recall the definition of $\mathfrak{K}_2$
introduced by D. Hadwin and J. Shen \cite{Don-Shen}, then we
introduce our new invariant $\mathfrak{K}_3$.

Let ${\cal M}_k({\Bbb C})$ be the $k\times k$ full matrix algebra
with entries in $\Bbb
 C$, and $\tau_k$ be the normalized
 trace on ${\cal M}_k({\Bbb C})$, i.e., $\tau_k=\frac{1}{k}Tr$,
 where $Tr$ is the usual trace on ${\cal M}_k({\Bbb C})$. Let
 ${\mathcal U}_k$ be the group of all unitary matrices in ${\cal M}_k(\Bbb
 C)$ and ${\cal M}_k({\Bbb C})^n$ denote the direct sum of $n$
 copies of ${\cal M}_k({\Bbb C})$. Define $\|\cdot\|_2$ on ${\cal M}_k(\Bbb
 C)^n$ by
 $$\|(A_1,\ldots,A_n)\|^2_2=\tau_k(A_1^*A_1)+\ldots+\tau_k(A_n^*A_n)$$for
 all $(A_1,\ldots,A_n)$ in ${\cal M}_k({\Bbb C})^n$.

 For every $\omega>0$, define the {\it $\omega$-ball}
 $Ball(B_1,\ldots,B_n;\omega)$ in ${\cal M}_k({\Bbb C})^n$ to be the subset of ${\cal M}_k(\Bbb
 C)^n$ consisting of all $(A_1,\ldots,A_n)$ in ${\cal M}_k(\Bbb
 C)^n$ such that $$\|(A_1,\ldots,A_n)-(B_1,\ldots,B_n)\|_2<\omega.$$
 Define the {\it $\omega$-orbit-ball} ${\mathcal
 U}(B_1,\ldots,B_n;\omega)$ in ${\cal M}_k({\Bbb C})^n$ to be the subset of ${\cal M}_k(\Bbb
 C)^n$ consisting of all $(A_1,\ldots,A_n)$ in ${\cal M}_k(\Bbb
 C)^n$ such that there exists some unitary matrix $W$ in ${\mathcal U}_k$ satisfying
  $$\|(A_1,\ldots,A_n)-(WB_1W^*,\ldots,WB_nW^*)\|_2<\omega.$$

Suppose $E\subseteq {\cal M}_k({\Bbb C})^n$, $\omega>0$. Define the
{\it covering number} $\nu_2(E,\omega)$ to be the minimal number of
$\omega$-balls that cover $E$ with the centers of these
$\omega$-balls in $E$; define the {\it $\omega$-orbit covering
number} $\nu(E,\omega)$ to be the minimal number of
$\omega$-orbit-balls that cover $E$ with the centers of these
$\omega$-orbit-balls in $E$.

Let $\mathcal M$ be a von Neumann algebra with a tracial state
$\tau$ and $x_1, x_2, \ldots, x_n$ be elements in ${\mathcal M}$.
For any $R,\varepsilon>0$, and positive integers $m$ and $k$, define
$\Gamma_R(x_1,\ldots,x_n;m,k,\varepsilon)$ to be the subset of
${\cal M}_k({\Bbb C})^n$ consisting of all $(A_1,\ldots,A_n)$ in
${\cal M}_k({\Bbb C})^n$ such that $\|A_j\|\leq R$ for $1\leq j\leq
n$, and
$$|\tau_k(A_{i_1}^{\eta_1}\cdots
A_{i_q}^{\eta_q})-\tau(x_{i_1}^{\eta_1}\cdots
x_{i_q}^{\eta_q})|<\varepsilon,$$ for all $1\leq i_1,\ldots,i_q\leq
n$, all $\eta_1,\ldots,\eta_q\in\{1,*\}$ and all $q$ with $1\leq
q\leq m$.

Define
$$\mathfrak{K}(x_1,\ldots,x_n;m,\varepsilon, \omega,R)=
\limsup_{k\rightarrow\infty}
\frac{\log(\nu(\Gamma_R(x_1,\ldots,x_n;m,k,\varepsilon),\omega))}{-k^2\log\omega}$$
$$\mathfrak{K}(x_1,\ldots,x_n;\omega,R)=
\inf_{m\in\Bbb N,\varepsilon>0}
\mathfrak{K}(x_1,\ldots,x_n;m,\varepsilon, \omega,R)$$
$$\mathfrak{K}(x_1,\ldots,x_n;\omega)=\sup_{R>0}\mathfrak{K}(x_1,\ldots,x_n;\omega,R)$$
$$\mathfrak{K}_2(x_1,\ldots,x_n)=\sup_{0< \omega< 1}\mathfrak{K}(x_1,\ldots,x_n;\omega)$$

D. Hadwin and J. Shen \cite{Don-Shen} also defined
$\mathfrak{K}_2(x_1,\ldots,x_n:y_1,\ldots,y_p)$ for all $x_1,
\ldots, x_n$, $y_1, \ldots, y_p$ in the von Neumann algebra $\cal M$
as follows. Let
\begin{eqnarray*}&&\Gamma_R(x_1,\ldots,x_n:y_1,\ldots,y_p;m,k,\varepsilon)\\
&=&\left\{(A_1,\ldots,A_n)\in {\cal M}_k({\Bbb C})^n:\ \mbox{there
exist}\ B_1,\ldots,B_p\ \mbox{in}\ {\cal M}_k({\Bbb C})\right.\\
&&\left.\mbox{such that}\
 (A_1,\ldots,A_n,B_1,\ldots,B_p)\in \Gamma_R(x_1,\ldots,x_n,y_1,\ldots,y_p;m,k,\varepsilon)
\right\},\end{eqnarray*}
$$\mathfrak{K}(x_1,\ldots,x_n:y_1,\ldots,y_p;m,\varepsilon,\omega,R)=
\limsup_{k\rightarrow\infty}
\frac{\log(\nu(\Gamma_R(x_1,\ldots,x_n;m,k,\varepsilon),\omega))}{-k^2\log\omega}$$
$$\mathfrak{K}(x_1,\ldots,x_n:y_1,\ldots,y_p;\omega,R)=
\inf_{m\in\Bbb N,\varepsilon>0}
\mathfrak{K}(x_1,\ldots,x_n:y_1,\ldots,y_p;m,\varepsilon,\omega,R)$$
$$\mathfrak{K}(x_1,\ldots,x_n:y_1,\ldots,y_p;\omega)=\sup_{R>0}\mathfrak{K}(x_1,\ldots,x_n;\omega,R)$$
$$\mathfrak{K}_2(x_1,\ldots,x_n:y_1,\ldots,y_p)=\sup_{0< \omega< 1}\mathfrak{K}(x_1,\ldots,x_n;\omega)$$

\begin{remark}\label{remark,x:y>x:y,z}From the definition, it is clear that

(1) $\mathfrak{K}_2(x_1,\ldots,x_n:y_1,\ldots,y_p)\geq
\mathfrak{K}_2(x_1,\ldots,x_n:y_1,\ldots,y_p,y_{p+1}),$

(2) if $\mathfrak{K}_2(x_1,\ldots,x_n:x_1,\ldots,x_{n+j})=0$ ($j\geq
0 $), then
$$\mathfrak{K}_2(x_1,\ldots,x_{n-1}:x_1,\ldots,x_{n+j})=0.$$
\end{remark}

Let $\infty\cdot 0=0$. For any subset $\mathcal G$ of $\mathcal M$,
define
$$\mathfrak{K}_3(x_1,\ldots,x_n:{\mathcal
G})=\inf\left\{\infty\cdot\mathfrak{K}_2(x_1,\ldots,x_n:A):\ A\
\mbox{is a finite subset of\ }{\mathcal G}\right\},$$
$$\mathfrak{K}_3({\mathcal
G})=\sup_{\begin{array}{l}E\subseteq {\mathcal G}\\
E\ \mbox{is finte}\end{array}}\inf_{\begin{array}{l}F\subseteq {\mathcal G}\\
F\ \mbox{is finte}\end{array}}\infty\cdot\mathfrak{K}_2(E:F).$$

 When
$\mathcal G$ is finite, it is not difficult to see that
$$\mathfrak{K}_3(x_1,\ldots,x_n:{\mathcal
G})=\infty\cdot\mathfrak{K}_2(x_1,\ldots,x_n:{\mathcal G})$$ and
$$\mathfrak{K}_3({\mathcal G})=\infty\cdot\mathfrak{K}_2({\mathcal
G}).$$

Note that the value of
$\mathfrak{K}_3(x_1,\ldots,x_n:y_1,\ldots,y_p)$ or
$\mathfrak{K}_3(x_1,\ldots,x_n)$ is always 0 or $\infty.$

\section{Key properties of $\mathfrak{K}_3$}
\begin{theorem}\label{equivalence}
If $\mathcal M$ is a von Neumann algebra with a tracial state
$\tau$, then the following are equivalent:

(1) $\mathfrak{K}_3({\mathcal M})=0$;

(2) if $x_1,\ldots,x_n\in\mathcal M$, then there exist
$y_1,\ldots,y_t\in\mathcal M$ such that
$\mathfrak{K}_2(x_1,\ldots,x_n:y_1,\ldots,y_t)=0$;

(3) for any generating set $\mathcal G$ of $\mathcal M$,
$\mathfrak{K}_3({\mathcal G})=0$;

(4) there exists a generating set $\mathcal G$ of $\mathcal M$ such
that $\mathfrak{K}_3({\mathcal G})=0$;

(5) if $\mathcal G$ is a generating set of $\mathcal M$, and $A_0$
is a finite subset of $\mathcal G$, then, for any finite subset $A$
with $A_0\subseteq A\subseteq \mathcal G$, there exists a finite
subset $B$ of $\mathcal G$ so that $\mathfrak{K}_2(A:B)=0$;

(6) there is an increasing directed family $\{{\mathcal
M}_\iota:\iota\in \Lambda\}$ of von Neumann subalgebras of $\mathcal
M$ such that

\ \ (a) each ${\mathcal M}_\iota$ is countably generated,

\ \ (b) $\mathfrak{K}_3({\mathcal M}_\iota)=0$,

\ \ (c) ${\mathcal M}=\cup_{\iota\in \Lambda}{\mathcal M}_{\iota}$.
\end{theorem}

Proof. It is clear that (1)$\Leftrightarrow$(2), (3)$\Rightarrow$(4)
and (3)$\Rightarrow$(5).

(4)$\Rightarrow$(2) Suppose $\cal G$ is a generating set of $\cal M$
and $\mathfrak{K}_3({\cal G})=0$. Let $\omega>0$ and $x_1, \ldots,
x_n$ be any elements in ${\cal M}$. Then there exist polynomials
$p_1, \ldots, p_n$ and elements $y_1,\ldots, y_s$ in $\cal G$ such
that
$$\|(x_1, \ldots, x_n)-(p_1(y_1,\ldots, y_s),\ldots, p_n(y_1,\ldots, y_s))\|\leq
\frac{\omega}{4}.$$

For any given $R>0$, if $(A_1, \ldots, A_s), (B_1, \ldots, B_s)$ in
${\cal M}_k({\Bbb C})^s$ and $\|A_j\|,\|B_j\|\leq R$ for all $1\leq
j\leq s$, then there exists a positive integer $N$ such that
\begin{eqnarray*}&&\|(p_1(A_1, \ldots, A_s), \ldots, p_n(A_1, \ldots,
A_s))-(p_1(B_1, \ldots, B_s),\ldots, p_n(B_1, \ldots, B_s))\|_2\\
& \leq & N\|(A_1, \ldots, A_s)-(B_1, \ldots,
B_s)\|_2.\end{eqnarray*}

Since $\mathfrak{K}_3({\cal G})=0$ and $y_1, \ldots, y_m$ are in
$\cal G$, there exist $y_1, \ldots, y_t$ ($t\geq s$) in $\cal G$
such that
$$\mathfrak{K}_2(y_1,\ldots, y_s:y_1, \ldots, y_t)=0.$$

Let $R> \max\{\|y_i\|, \|x_j\|: 1\leq i\leq n, 1\leq j\leq t\}$.
Note that if $(A_1,\ldots, A_n,B_1, \ldots, B_t)\in \Gamma_R(x_1,
\ldots, x_n,y_1, \ldots, y_t;m,k,\varepsilon)$, then, for
sufficiently small $\varepsilon$ and sufficiently large $m$, we have
$$\|(A_1, \ldots, A_n)-(p_1(B_1,\ldots, B_s),\ldots, p_n(B_1,\ldots, B_s))\|_2<\frac{\omega}{4}$$
and $$(B_1, \ldots, B_s)\in\Gamma(y_1, \ldots, y_s:y_1, \ldots,
y_t;m,k,\varepsilon).$$ It follows that there is a set $\Lambda$ and
a subset $\{(B_1^{\lambda}, \ldots,
B_s^{\lambda}):\lambda\in\Lambda\}$ of $\Gamma_R(y_1, \ldots,
y_s:y_1, \ldots, y_t;m,k,\varepsilon)$ with
$$\mbox{card}(\Lambda)\leq \nu(\Gamma(y_1, \ldots, y_s:y_1, \ldots,
y_t;m,k,\varepsilon),\frac{\omega}{4N}).$$ That means, for every
$(A_1, \ldots, A_n,B_1, \ldots, B_t)\in \Gamma_R (x_1, \ldots,
x_n,y_1, \ldots, y_t;m,k,\varepsilon)$, there is a
$\lambda\in\Lambda$ and a unitary $k\times k$ matrix $U$ such that
$$\|(B_1, \ldots, B_s)-U^*(B_1^{\lambda}, \ldots,
B_s^{\lambda})U\|_2\leq \frac{\omega}{4N}.$$ That gives
$$\|(A_1, \ldots A_n)-U^*(p_1(B_1^{\lambda}, \ldots, B_s^{\lambda}),\ldots, p_n(B_1^{\lambda}, \ldots, B_s^{\lambda}))U\|_2<\frac{\omega}{2}.$$

It follows that, if $\varepsilon$ is sufficient small and $m$ is
sufficient large, then for any $k\in\Bbb N$,
\begin{eqnarray*}&&\frac{\nu(\Gamma_R(x_1, \ldots, x_n:y_1, \ldots, y_t;m,k,\varepsilon),\omega)}{-k^2\log \omega}\\
&\leq&
\left(\frac{\log\omega}{\log(\omega/4N)}\right)\frac{\nu(\Gamma(y_1,
\ldots, y_s:y_1, \ldots,
y_t;m,k,\varepsilon),\omega/4N)}{-k^2\log(\omega/4N)}.\end{eqnarray*}
If we take $k\rightarrow \infty$ and take the infimum over $m$ and
$\varepsilon$, we get
\begin{eqnarray*}&&\mathfrak{K}_2(x_1, \ldots, x_n:y_1, \ldots, y_t;\omega)\\
&\leq&
\left(\frac{\log\omega}{\log(\omega/4N)}\right)\mathfrak{K}_2(y_1,
\ldots, y_s:y_1, \ldots, y_t)\\
&=&0.\end{eqnarray*}

Hence, we have $\mathfrak{K}_2(x_1, \ldots, x_n: y_1, \ldots,
y_t)=0$.

(2)$\Rightarrow$(3)\ \ Suppose (2) is true, and suppose
$x_1,\ldots,x_n$ are elements of some generating set $\mathcal G$ of
$\mathcal M$. Then there exist $y_1,\ldots,y_t$ in $\mathcal M$ such
that $\mathfrak{K}_2(x_1,\ldots,x_n:y_1,\ldots,y_t)=0$.

Suppose $\varepsilon_0>0$, $m_0\in \Bbb N$, and $0< \omega< 1$. We
can choose $w_1,\ldots,w_s$ in $\mathcal G$ and polynomials
$p_1,\ldots,p_t$ so that each $\|y_j-p_j(w_1,\ldots, w_s)\|_2$
($1\leq j\leq t$) is small enough to make
$$\left|\tau(q(x_1,\ldots,x_n,y_1,\ldots,y_t))-
\tau(q(x_1,\ldots,x_n,p_1(w_1,\ldots,w_s),\ldots,p_t(w_1,\ldots,w_s)))\right|<\frac{\varepsilon_0}{4},$$
for every monomial $q$ with length at most $m_0$.

When $m$ is sufficient large, $\varepsilon$ is sufficient small,
if
$$(A_1,\ldots,A_n,B_1,\ldots,B_s)\in\Gamma_R(x_1,\ldots,x_n,w_1,\ldots,w_s;m,k,\varepsilon),$$
 then $$(A_1,\ldots,A_n,p_1(B_1,\ldots,B_s),\ldots,p_t(B_1,\ldots,B_s))\in
 \Gamma_R(x_1,\ldots,x_n,y_1,\ldots,y_t;m_0,k,\varepsilon_0).$$
Hence
$$\Gamma_R(x_1,\ldots,x_n:w_1,\ldots,w_s;m,k,\varepsilon)\subseteq
\Gamma_R(x_1,\ldots,x_n:y_1,\ldots,y_t;m_0,k,\varepsilon_0).$$ Since
$$\nu(\Gamma_R(x_1,\ldots,x_n:w_1,\ldots,w_s;m,k,\varepsilon),\omega)\leq
2\nu(\Gamma_R(x_1,\ldots,x_n:y_1,\ldots,y_t;m_0,k,\varepsilon_0),\omega),$$
we have,
$$\mathfrak{K}(x_1,\ldots,x_n:w_1,\ldots,w_s;\omega)\leq
2\mathfrak{K}(x_1,\ldots,x_n:y_1,\ldots,y_t;m_0,\varepsilon_0,\omega).$$
Then we get
$$\mathfrak{K}(x_1,\ldots,x_n:w_1,\ldots,w_s;\omega)\leq
2\mathfrak{K}(x_1,\ldots,x_n:y_1,\ldots,y_t;\omega).$$ Therefore
$$\mathfrak{K}_2(x_1,\ldots,x_n:w_1,\ldots,w_s)\leq 2\mathfrak{K}_2(x_1,\ldots,x_n:y_1,\ldots,y_t)=0.
$$Thus
$\mathfrak{K}_2(x_1,\ldots,x_n:w_1,\ldots,w_s)=0$. From the
definition, $\mathfrak{K}_3({\mathcal G})=0$.

(5)$\Rightarrow$(6)\ \ Suppose $A_0\subseteq
\{x_1,\ldots,x_n\}\subseteq \mathcal G$. From (4), there exists a
family $\{B_0, B_1, \ldots\}$ of finite subsets of $\mathcal G$ such
that $$\mathfrak{K}_2(x_1, \ldots, x_n: B_0)=0,$$ and for any
positive integer $n$, $$\mathfrak{K}_2(x_1, \ldots, x_n, B_0,
\ldots, B_{n-1}: B_n)=0.$$ Let $\cal G$ be the set $\{x_1,\ldots,
x_n\}\cup\cup_{n=0}^{\infty}B_n$ and $\cal N$ be the von Neumann
subalgebra generated by $\cal G$. Then $\cal N$ is countably
generated.

Let $A\subseteq \cal G$ be a finite subset. Then there exists a
positive integer $m$, so that $A\subseteq \{x_1, \ldots, x_n\}\cup
B_0\cup\cdots \cup B_m$. Since
$$\mathfrak{K}_2(x_1, \ldots, x_n, B_0, \ldots, B_{m}: B_{m+1})=0,$$
by Remark \ref{remark,x:y>x:y,z}, we have
$$\mathfrak{K}_2(A: B_n)=0.$$It follows that
$\mathfrak{K}_3\left({\cal G}\right)=0.$ Therefore
$\mathfrak{K}_3({\cal N})=0$ by the equivalence of (1) and (4).

It is not difficult to see that the union of all such $\mathcal N$'s
is $\mathcal M$.

(6)$\Rightarrow$(2)\ \ Suppose $x_1,\ldots,x_n$ are elements of
$\mathcal M$. From (5), there exists $\{\iota_1,\ldots,
\iota_n\}\subseteq \Lambda$ such that $x_1\in {\mathcal
M}_{\iota_1},\ldots, x_n\in{\mathcal M}_{\iota_n}$. Since
$\{{\mathcal M}_\iota:\iota\in \Lambda\}$ is an increasing directed
family, there exists $\iota\in\Lambda$, such that ${\mathcal
M}_{\iota_1},\ldots,{\mathcal M}_{\iota_n}\subseteq {\cal
M}_{\iota}$ and $\mathfrak{K}_3({\cal M}_{\iota})=0$. Therefore
there exist $y_1,\ldots,y_m\in {\mathcal M}_{\iota}\subseteq
\mathcal M$ such that $\mathfrak{K}_2(x_1,\ldots, x_n:y_1,\ldots,
y_m)=0$.\hfill $\Box$

\begin{remark}If $\cal M$ is finite generated, then $\mathfrak{K}_3({\cal M})=0$ is equivalent to $\mathfrak{K}_2({\cal M})=0$.
\end{remark}

\begin{corollary}\label{corollary, k3 is independent on generatoer}Suppose $\mathcal M$
 is a von Neumann algebra with a tracial state $\tau$, $\mathcal G$ is
a generating set of $\mathcal M$. Then $\mathfrak{K}_3({\mathcal M
})=\mathfrak{K}_3({\mathcal G})$.
\end{corollary}

\begin{corollary}\label{corollary, increasingly directed}Suppose $\{{\mathcal M}_{\iota}\}_{\iota\in\Lambda}$
is an increasingly directed family of von Neumann algebras. Then
$\mathfrak{K}_3(\cup_{\iota}{\mathcal M}_{\iota})\leq
\liminf_{\iota}\mathfrak{K}_3({\mathcal M}_{\iota}).$
\end{corollary}

\begin{remark}To see that Corollary \ref{corollary, increasingly directed} gives the best estimate, note that ${\mathcal L}_{{\Bbb F}_2}\otimes
{\mathcal R}=\cup_n\left( {\mathcal L}_{{\Bbb F}_2}\otimes{\mathcal
M}_{2^n}({\Bbb C})\right)$ and $\mathfrak{K}_3({\mathcal L}_{{\Bbb
F}_2}\otimes {\mathcal R})=0$, but $\mathfrak{K}_3({\mathcal
L}_{{\Bbb F}_2}\otimes {\mathcal M}_{2^n}(\Bbb C))=\infty$ for every
$n$.
\end{remark}

\vspace{0.2cm}

To prove Theorem \ref{theorem, intersection contains Haar unitary},
we need the following lemmas.
\begin{lemma}\label{add generator is the same} Let $\mathcal M$ be a von Neumann algebra with a tracial state $\tau$. Suppose
$x_1$, $\ldots$, $x_n$, $y_1$, $\ldots$, $y_p$, $w_1,\ldots,w_t$ are
elements of $\mathcal M$ and $x_1,\ldots,x_n\in
W^*(y_1,\ldots,y_p)$. Then, for $\omega>0$,
$$\mathfrak{K}(y_1,\ldots,y_p:w_1,\ldots,w_t;\omega)=\mathfrak{K}(y_1,\ldots,y_p:x_1,\ldots,x_n,w_1,\ldots,w_t;\omega).$$
\end{lemma}

Proof. It is not hard to get the ``$\geq$" part.

Assume $\varepsilon_0>0$, $m_0\in\Bbb N$, $R>1$. Since
$x_1,\ldots,x_n\in W^*(y_1,\ldots,y_p)$, there exist $m_1\in \Bbb N$
and $\varepsilon_1>0$ and a family of noncommutative polynomials
$q_1,\ldots, q_n$ such that
\newline $\|\left(q_1(y_1,\ldots,y_p),\ldots,
q_n(y_1,\ldots,y_p)\right)-(x_1,\ldots,x_n)\|_2$ is so small that
for any $m\geq m_1$ and $0<\varepsilon\leq \varepsilon_1$, we have,
for any $k\in\Bbb N$,
\begin{eqnarray*}&&\{(A_1,\ldots,A_p,q_1(A_1,\ldots,A_p),\ldots, q_n(A_1,\ldots,A_p),C_1,\ldots,C_t):\\
&&(A_1,\ldots,A_p,C_1,\ldots,C_t)\in\Gamma_R(y_1,\ldots,y_p,w_1,\ldots,w_t;m,k,\varepsilon)\}\\
&\subseteq&
\Gamma_R(y_1,\ldots,y_p,x_1,\ldots,x_n,w_1,\ldots,w_t;m_0,k,\varepsilon_0),
\end{eqnarray*}
which implies
$$\Gamma_R(y_1,\ldots,y_p:w_1,\ldots,w_t;m,k,\varepsilon)\subseteq
\Gamma_R(y_1,\ldots,y_p:x_1,\ldots,x_n,w_1,\ldots,w_t;m_0,k,\varepsilon_0)
.$$ Therefore
$\mathfrak{K}(y_1,\ldots,y_p:w_1,\ldots,w_t;\omega)\leq\mathfrak{K}(y_1,\ldots,y_p:x_1,\ldots,x_n,w_1,\ldots,w_t;\omega)$.\hfill
$\Box$

\vspace{0.2cm}

The following lemma is a slight extension of Theorem 1 in
\cite{Don-Shen}; the proofs are similar.
\begin{lemma}\label{change generator is the same} Let
$x_1,\ldots,x_n,y_1,\ldots,y_p,w_1,\ldots,w_t$ be elements in a von
Neumann algebra $\mathcal M$ with a tracial state $\tau$, and
$W^*(x_1,\ldots,x_n)=W^*(y_1,\ldots,y_p)$. Then
$$\mathfrak{K}_3(x_1,\ldots,x_n:w_1,\ldots,w_t)=\mathfrak{K}_3(y_1,\ldots,y_p:w_1,\ldots,w_t).$$
\end{lemma}

Proof. If $\mathfrak{K}_3(x_1,\ldots,x_n:w_1,\ldots,w_t)$ and
$\mathfrak{K}_3(y_1,\ldots,y_p:w_1,\ldots,w_t)$ are both infinity,
then they are equal. If one of them is zero, say it
$\mathfrak{K}_3(x_1,\ldots,x_n:w_1,\ldots,w_t)$, then we need to
prove that $\mathfrak{K}_3(y_1,\ldots,y_p:w_1,\ldots,w_t)=0$.

For every $0<\omega<1$, there exists a family of noncommutative
polynomials $q_1,\ldots, q_p$, such that
$$\|(y_1,\ldots,y_p)-\left(q_1(x_1,\ldots,x_n),\ldots,q_p(x_1,\ldots,x_n\right))\|_2\leq\frac{\omega}{4}.$$
For such a family of polynomials $q_1,\ldots,q_p$ and every $R>0$,
there always exists a constant $D\geq 1$, depending only on
$q_1,\ldots,q_p$ and $R$, such that
\begin{eqnarray*}&&\|\left(q_1(A_1,\ldots,A_n),\ldots,q_p(A_1,\ldots,A_n)\right)-\left(q_1(B_1,\ldots,B_n),\ldots,q_p(B_1,\ldots,B_n)\right)\|_2\\
&\leq& D\|(A_1,\ldots,A_n)-(B_1,\ldots,B_n)\|_2,\end{eqnarray*} for
all $A_1,\ldots,A_n,B_1,\ldots,B_n$ in ${\mathcal M}_k({\Bbb C})$
satisfying $\|A_j\|, \|B_j\|\leq R$, for $1\leq j\leq n$.

For $R>1$, $m$ and $k$ sufficiently large, $\varepsilon$
sufficiently small, if
$$(B_1,\ldots,B_p,A_1,\ldots,A_n)\in
\Gamma_R(y_1,\ldots,y_p,
x_1,\ldots,x_n:w_1,\ldots,w_t;m,k,\varepsilon),$$ then
$$\|(B_1,\ldots,B_p)-\left(q_1(A_1,\ldots,A_n),\ldots,q_p(A_1,\ldots,A_n)\right)\|_2\leq \frac{\omega}{4}.$$
It is clear that $(A_1,\ldots,A_n)\in \Gamma_R(
x_1,\ldots,x_n:w_1,\ldots,w_t;m,k,\varepsilon)$.

There exists a set $\{\mathcal
U(A_1^{\lambda},\ldots,A_n^{\lambda};\frac{\omega}{4D})\}_{\lambda\in\Lambda_k}$
of $\frac{\omega}{4D}$-orbit-balls that cover $\Gamma_R(
x_1,\ldots,x_n:w_1,\ldots,w_t;\\ m,k,\varepsilon)$ with the
cardinality of $\Lambda_k$ satisfying $|\Lambda_k|=\nu(\Gamma_R(
x_1,\ldots,x_n:w_1,\ldots,w_t;m,k,\varepsilon),\frac{\omega}{4D})$.
Thus there exists some $\lambda\in\Lambda_k$ and $U\in {\mathcal
U}_k$ such that
$$\|(A_1,\ldots,A_n)-U(A_1^{\lambda},\ldots,A_n^{\lambda})U^*\|_2\leq
\frac{\omega}{4D}.$$ It follows that
\begin{eqnarray*}&&\|(B_1,\ldots,B_p)-U\left(q_1(A_1^{\lambda},\ldots,A_n^{\lambda}),\ldots,q_p(A_1^{\lambda},\ldots,A_n^{\lambda})\right)U^*\|_2\\
&=&\|(B_1,\ldots,B_p)-\left(q_1(U(A_1^{\lambda},\ldots,A_n^{\lambda})U^*),\ldots,q_p(U(A_1^{\lambda},\ldots,A_n^{\lambda})U^*)\right)\|_2\\
&\leq & \frac{\omega}{2}.\end{eqnarray*} That is,
$$(B_1,\ldots,B_p)\in{\mathcal
U}\left(q_1(A_1^{\lambda},\ldots,A_n^{\lambda}),\ldots,q_p(A_1^{\lambda},\ldots,A_n^{\lambda});\omega\right).$$
Hence, we get \begin{eqnarray*} 0&\leq &\mathfrak{K}(y_1,\ldots,y_p:x_1,\ldots,x_n,w_1,\ldots,w_t;\omega, R)\\
&\leq & \inf_{m\in\Bbb N,\varepsilon>0}\limsup_{k\rightarrow \infty}\frac{\log(|\Lambda_k|)}{-k^2\log\omega}\\
&=&\inf_{m\in\Bbb N,\varepsilon>0}\limsup_{k\rightarrow
\infty}\frac{\log(\nu(\Gamma_R(
x_1,\ldots,x_n:w_1,\ldots,w_t;m,k,\varepsilon),\frac{\omega}{4D}))}{-k^2\log\omega}\\
&=&0.
\end{eqnarray*}
Therefore
$\mathfrak{K}_2(y_1,\ldots,y_p:x_1,\ldots,x_n,w_1,\ldots,w_t)=0$.
From Lemma \ref{add generator is the same}, we get
$\mathfrak{K}_2(y_1,\ldots,y_p:w_1,\ldots,w_t)=0$. So
$\mathfrak{K}_3(y_1,\ldots,y_p:w_1,\ldots,w_t)=0$.\hfill $\Box$

\begin{definition}A unitary matrix $A$ in ${\mathcal M}_k({\Bbb C})$
is called a Haar unitary matrix if the eigenvalues of $A$ are the
$k$-th roots of unity; equivalently, if $\tau_k(A^i)=0$ for $1\leq
i<k$ and $\tau_k(A^k)=1$.
\end{definition}

\begin{lemma}[\cite{Vo}]\label{Vo} Let $V_1,V_2$ be two Haar unitary matrices in ${\mathcal M}_k(\Bbb
C)$. For every $\delta>0$, let $$\Omega(V_1,V_2;\delta)=\{U\in
{\mathcal U}_k: \|UV_1-V_2U\|_2\leq \delta\}.$$ Then, for every
$0<\delta<r$, $\nu_2(\Omega(V_1,V_2;\delta),\frac{4\delta}{r})\leq
(\frac{3r}{2\delta})^{4rk^2}$.
\end{lemma}

\begin{definition}Suppose $\mathcal M$ is a von Neumann algebra
 with a tracial state $\tau$. Then a unitary $u$ in $\mathcal M$
 is called a Haar unitary if $\tau(u^m)=0$ when $m\neq 0$. In
 addition, $\cal M$ is called diffuse if $\cal M$ contains a Haar
 unitary.
\end{definition}

The following lemma is a slight extension of Theorem 6 in
\cite{Don-Shen}, and the proofs are similar.
\begin{lemma}\label{haar unitary} Let $x_1,\ldots,x_n, y_1,\ldots,y_p, v_1,\ldots,v_s, w_1,\ldots,w_t$
 be elements in a von Neumann algebra $\mathcal M$ with a
tracial state $\tau$. If $W^*(x_1,\ldots,x_n)\cap
W^*(y_1,\ldots,y_p)$ is diffuse, then
\begin{eqnarray*}&&\mathfrak{K}_3(x_1,\ldots,x_n,y_1,\ldots,y_p:v_1,\ldots,v_s,w_1,\ldots,w_t)\\
&\leq & \mathfrak{K}_3(x_1,\ldots,x_n:v_1,\ldots,v_s)+
\mathfrak{K}_3(y_1,\ldots,y_p:w_1,\ldots,w_t). \end{eqnarray*}
\end{lemma}

Proof. If one of $\mathfrak{K}_3(x_1,\ldots,x_n:v_1,\ldots,v_s)$ and
$\mathfrak{K}_3(y_1,\ldots,y_p:w_1,\ldots,w_t)$ is infinity, then we
are done.

Now suppose
$\mathfrak{K}_3(x_1,\ldots,x_n:v_1,\ldots,v_s)=\mathfrak{K}_3(y_1,\ldots,y_p:w_1,\ldots,w_t)=0$.
Since $W^*(x_1,\ldots,x_n)\cap W^*(y_1,\ldots,y_p)$ is diffuse, we
can find a Haar unitary $u$ in $W^*(x_1,\ldots,x_n)\cap
W^*(y_1,\ldots,y_p)$.

For $R>1+\max_{1\leq i\leq n, 1\leq j\leq p}\{\|x_i\|,\|y_j\|\}$,
$0<\omega<\frac{1}{2n}$, $0<r<1$ and $\varepsilon>0$, $m,k\in\Bbb
N$. Suppose
$$(A_1,\ldots,A_n, B_1,\ldots,B_p,U)\in\Gamma_R(x_1,\ldots,x_n,
y_1,\ldots,y_p,u: v_1,\ldots,v_s, w_1,\ldots,w_t;m,k,\varepsilon).$$
Then
$$(A_1,\ldots,A_n,U)\in\Gamma_R(x_1,\ldots,x_n,u:v_1,\ldots,v_s,w_1,\ldots,w_t;m,k,\varepsilon)$$
and
$$(B_1,\ldots,B_p,U)\in\Gamma_R(y_1,\ldots,y_p,u:v_1,\ldots,v_s,w_1,\ldots,w_t;m,k,\varepsilon).$$

Let $\{{\mathcal
U(A_1^{\lambda},\ldots,A_n^{\lambda},U^{\lambda})};\frac{r\omega}{24R}\}_{\lambda\in\Lambda_k}$
be a set of $\frac{r\omega}{24R}$-orbit-balls that cover
$\Gamma_R(x_1,\ldots,x_n,u: v_1,\ldots,v_s,\\
w_1,\ldots,w_t;m,k,\varepsilon)$ with the cardinality of $\Lambda_k$
satisfying
$$|\Lambda_k|=\nu(\Gamma_R(x_1,\ldots,x_n,u:
v_1,\ldots,v_s,w_1,\ldots,w_t;m,k,\varepsilon);\frac{r\omega}{24R}).$$
Also let $\{{\mathcal
U(B_1^{\sigma},\ldots,B_p^{\sigma},U^{\sigma})};\frac{r\omega}{24R}\}_{\lambda\in\Sigma_k}$
be a set of $\frac{r\omega}{24R}$-orbit-balls that cover
$\Gamma_R(y_1,\ldots,y_p,u:v_1,\ldots,v_s,\\
w_1,\ldots,w_t;m,k,\varepsilon)$ with the cardinality of $\Sigma_k$
satisfying
$$|\Sigma_k|=\nu(\Gamma_R(y_1,\ldots,y_p,u:v_1,\ldots,v_s,w_1,\ldots,w_t;m,k,\varepsilon);\frac{r\omega}{24R}).$$
When $m$ is sufficiently large and $\varepsilon$ is sufficiently
small, by Theorem 2.1 in \cite{Don1}, we can assume that all
$U^{\lambda}, U^{\sigma}$ to be Haar unitary matrices in ${\cal
M}_k(\Bbb C)$.

For any $$(A_1,\ldots,A_n,
B_1,\ldots,B_p,U)\in\Gamma_R(x_1,\ldots,x_n, y_1,\ldots,y_p,u:
v_1,\ldots,v_s, w_1,\ldots,w_t;m,k,\varepsilon),$$ there exist some
$\lambda\in \Lambda_k$, $\sigma\in \Sigma_k$ and
$W_1,W_2\in{\mathcal U }_k$ such that
$$\|(A_1,\ldots,A_n,U)-W_1(A_1^{\lambda},\ldots,A_n^{\lambda},U^{\lambda})W_1^*\|_2\leq \frac{r\omega}{24R},$$
$$\|(B_1,\ldots,B_p,U)-W_2(B_1^{\sigma},\ldots,B_p^{\sigma},U^{\sigma})W_2^*\|_2\leq \frac{r\omega}{24R}.$$
Therefore
$$\|W_1U^{\lambda}W_1^*-W_2U^{\sigma}W_2^*\|_2=\|W_2^*W_1U^{\lambda}-U^{\sigma}W_2^*W_1\|_2<\frac{r\omega}{12R}.$$

From Lemma \ref{Vo}, there exists a set
$\{Ball(U_{\lambda,\sigma,\gamma},\frac{\omega}{3R})\}_{\gamma\in
\Delta_k}$ in ${\mathcal U}_k$ which cover
$\Omega(U^{\lambda},U^{\sigma};\frac{r\omega}{12R})$ with
cardinality $|\Delta_k|\leq (\frac{18R}{\omega})^{4rk^2}$. This
implies
\begin{eqnarray*}&&\|(A_1,\ldots,A_n,B_1,\ldots,B_p,U)\\
&-&
(W_2U_{\lambda,\sigma,\gamma}A_1^{\lambda}U_{\lambda,\sigma,\gamma}^*W_2^*,
\ldots,W_2U_{\lambda,\sigma,\gamma}A_n^{\lambda}U_{\lambda,\sigma,\gamma}^*W_2^*,
W_2B_1^{\sigma}W_2^*,\ldots,W_2B_p^{\sigma}W_2^*,W_2U^{\sigma}W_2^*)\|_2\\
&\leq& n\omega.\end{eqnarray*}
 Then we get
\begin{eqnarray*}&&\mathfrak{K}(x_1,\ldots,x_n,
y_1,\ldots,y_p,u: v_1,\ldots,v_s, w_1,\ldots,w_t;2n\omega,R)\\
&\leq&\inf_{m\in\Bbb N,\varepsilon>0}
\limsup_{k\rightarrow \infty}\frac{\log(|\Lambda_k||\Sigma_k||\Delta_k|)}{-k^2\log(2n\omega)}\\
&\leq &4r\frac{\log (18R)-\log \omega}{-\log(2n\omega)}.
\end{eqnarray*}

Because $r$ is an arbitrarily small positive number, we have
$$\mathfrak{K}_3(x_1,\ldots,x_n,y_1,\ldots,y_p,u:v_1,\ldots,v_s,w_1,\ldots,w_t)=0.$$
Note that
$W^*(x_1,\ldots,x_n,y_1,\ldots,y_p,u)=W^*(x_1,\ldots,x_n,y_1,\ldots,y_p)$,
by Lemma \ref{change generator is the same}, we have
$$\mathfrak{K}_3(x_1,\ldots,x_n,y_1,\ldots,y_p:v_1,\ldots,v_s,w_1,\ldots,w_t)=0.$$\hfill $\Box$

\vspace{0.2cm}

Now we are ready to prove the following theorem.

\begin{theorem}\label{theorem, intersection contains Haar unitary}Suppose $\mathcal M$ is a von Neumann algebra with a
tracial state $\tau$, ${\mathcal N}_1$ and ${\mathcal N}_2$ are von
Neumann subalgebras of $\mathcal M$. If ${\mathcal N}_1\cap{\mathcal
N}_2$ is diffuse, then
$$\mathfrak{K}_3(({\cal N}_1\cup {\cal N}_2)'')\leq
\mathfrak{K}_3({\mathcal N}_1)+\mathfrak{K}_3({\mathcal N}_2).$$
\end{theorem}

Proof. If one of $\mathfrak{K}_3({\mathcal N}_1)$ and
$\mathfrak{K}_3({\mathcal N}_2)$ is infinity, then we are done.

Now suppose $\mathfrak{K}_3({\mathcal N}_1)=\mathfrak{K}_3({\mathcal
N}_2)=0$ and $u$ is a Haar unitary in ${\mathcal N}_1\cap{\mathcal
N}_2$.  Let ${\cal G}={\mathcal N}_1\cup{\mathcal N}_2$ and
$A_0=\{u\}$. Then $\cal G$ is a generating set of $\mathcal M$.
Suppose $A_0\subseteq A\subseteq {\cal G}$ and $A$ is finite, write
 $A=\{x_1,\ldots,x_n,u,y_1,\ldots,y_p\}$ with $x_1,\ldots,x_n\in{\mathcal
N}_1$ and $y_1,\ldots,y_p\in{\mathcal N}_2$. Since
$\mathfrak{K}_3({\mathcal N}_1)=\mathfrak{K}_3({\mathcal N}_2)=0$,
there exist $v_1,\ldots,v_s\in {\mathcal N}_1$,
$w_1,\ldots,w_t\in{\mathcal N}_2$ such that
$\mathfrak{K}_2(x_1,\ldots,x_n,u:v_1,\ldots,v_s)=0$ and
$\mathfrak{K}_2(y_1,\ldots,y_p,u:w_1,\ldots,w_t)=0.$

 Because $u\in
W^*(x_1,\ldots,x_n,u)\cap W^*(y_1,\ldots,y_p,u)$, then from Lemma
\ref{haar unitary}, we know that
\begin{eqnarray*}&&\mathfrak{K}_2(A:v_1,\ldots,v_s,w_1,\ldots,w_t)\\
&=&\mathfrak{K}_2(x_1,\ldots,x_n,u,y_1,\ldots,y_p:v_1,\ldots,v_s,w_1,\ldots,w_t)\\
&=&0.\end{eqnarray*}

Therefore, by Theorem \ref{equivalence}, $\mathfrak{K}_3(({\cal
N}_1\cup {\cal N}_2)'')=0$.\hfill $\Box$

\begin{lemma}\label{don} (\cite{Don}, Lemma 17) Suppose $\mathcal M$ is a von Neumann algebra with a
tracial state $\tau$, $x$ is a normal element in $\mathcal M$ such
that $x$ has no eigenvalues. Then there is a selfadjoint element $y$
with the uniform distribution on $[0,1]$ such that $W^*(x)=W^*(y)$.
\end{lemma}

\begin{lemma}\label{don, lemma 18} (\cite{Don}, Lemma 18) Suppose $n,n_1,p\in\Bbb
N$, $1\leq n_1\leq n$ and $p\geq 2$. Suppose $A$ is a diagonal
matrix whose diagonal entries are $\frac{1}{n_1}, \frac{2}{n_1},
\cdots, \frac{n_1}{n_1}, -1, \cdots, -1$, and $B$ is any selfadjoint
$n\times n$ matrix with $0\leq B\leq 1$. Suppose $0\leq
\varepsilon\leq \frac{1}{4}$, and
$\frac{1}{n_1}<4\varepsilon^{p-1}$. Let
$$\Sigma(A,B,\varepsilon^p)=\{W\in {\cal M}_n: \|W\|\leq 1, \|AW-WB_2\|_2<\varepsilon^p\}.$$
Then $\nu_2(\Sigma(A,B,\varepsilon^p),\varepsilon)\leq
(\frac{6}{\varepsilon})^{16n^2\varepsilon^{p-1}}$.
\end{lemma}

\begin{theorem}\label{theorem, without comon eigenvalues} Let $\mathcal M$ be a von Neumann algebra with a
tracial state $\tau$. Suppose $\mathcal N$ is a von Neumann
subalgebra of $\mathcal M$ and $y$ is an element in $\mathcal M$. If
$a$, $b$ are normal operators in $\mathcal N$ such that $a$ and $b$
have no common eigenvalues and $ay=yb\neq 0$, then
$\mathfrak{K}_3\left(W^*({\cal N}\cup
\{y\})\right)\leq\mathfrak{K}_3({\mathcal N})$.
\end{theorem}

Proof. There is no loss in assuming that $y, a, b$ have norm at most
1.

If $\mathfrak{K}_3({\mathcal N})$ is infinity, then we are done. So
we can assume $\mathfrak{K}_3({\mathcal N})=0$.

Since $ay=yb$, by the Putnam-Fuglede theorem, we have $a^*y=yb^*$.
Then for any polynomial $p$, $p(a,a^*)a=ap(b,b^*)$.
% $b_1b_^*a=ab_2b_^*$
 Therefore for every bounded Borel
function $\varphi:{\Bbb C}\rightarrow{\Bbb C}$,
$\varphi(a)y=y\varphi(b)$. If $\lambda\in\sigma_p(a)$, then the
spectral projection $\chi_{\{\lambda\}}(b)$ is 0. If $a$ is
diagonal, then
$y=Iy=\sum_{\lambda\in\sigma_p(a)}\chi_{\{\lambda\}}(a)y=y\sum_{\lambda\in\sigma_p(a)}\chi_{\{\lambda\}}(b)=0.$
So $a$ can  not be diagonal, therefore we can write $a=d\oplus c$,
where $d$ is diagonal and $c$ has no eigenvalues. From Lemma
\ref{don}, there is a self-adjoint element $c_0$ with the uniform
distribution such that $W^*(c)=W^*(c_0)$, thus there is some Borel
function $\psi$ such that $c_0=\psi(c)$ and $\psi(d)=-1$, $0\leq
\psi(b)\leq 1$. Since $\psi(a)y=(-1\oplus c_0)y=y\psi(b)$, we can
replace $-1\oplus c_0$ with $a$, $\psi(b)$ with $b$. Hence we can
assume $b$ is selfadjoint and $a=-1\oplus c_0$ with
$\tau(a^n)=(-1)^n(1-\alpha)+\alpha\int_0^1t^ndt$, where
$\alpha=1-\tau(1\oplus 0).$

For each $n\in\Bbb N$, define the diagonal matrix $A_n$ with
eigenvalues
$$\frac{1}{[n\alpha]},\frac{2}{[n\alpha]},\ldots,\frac{[n\alpha]}{[n\alpha]},-1,-1,\ldots,-1,$$
where $[n\alpha]$ denotes the greatest integer function of
$n\alpha$. It is not hard to show that
$\tau_n(A_n^m)=(-1)^m\frac{n-[n\alpha]}{n}+
\frac{[n\alpha]}{m}\left(\frac{1}{[n\alpha]}\sum_{s=1}^{[n\alpha]}(\frac{s}{[n\alpha]})^m\right)$.
When $n\rightarrow\infty$,
$\frac{1}{[n\alpha]}\sum_{s=1}^{[n\alpha]}(\frac{s}{[n\alpha]})^m$
is a Riemann sum converging to $\int_0^1t^mdt.$ Hence
$\tau_n(f(A_n))\rightarrow \tau(f(a))$ as $n\rightarrow \infty$.
Choose matrix $B_n\in{\mathcal M}_n({\Bbb C})$ such that $0\leq
B_n\leq 1$ and $B_n$ converges in distribution to $b$ as
$n\rightarrow \infty$.

Let $x_1,\ldots,x_n$ be elements in $\mathcal N$. Then there exist
$y_1,\ldots,y_p$ in $\mathcal N$ such that
$$\mathfrak{K}_2(x_1,\ldots,x_n:y_1,\ldots,y_p)=0.$$ For any
$0<\omega<1$, $0<r<1$, $R>1$, $m\in\Bbb N$, $\varepsilon>0$ and
$k\in\Bbb N$, there exists a set $\{{\mathcal
U}(T_1^{\lambda},\ldots,T_n^{\lambda},A^{\lambda},B^{\lambda};\\\frac{r\omega}{64})\}_{\lambda\in\Lambda_k}$
 of $\frac{r\omega}{64}$-orbit-balls in ${\mathcal M}_k(\Bbb
 C)^{k+2}$ that cover
 $\Gamma_R(x_1,\ldots,x_n,a,b:y_1,\ldots,y_p;m,k,\varepsilon)$ with the
 cardinality of $\Lambda_k$ satisfying $|\Lambda_k|=\nu(\Gamma_R(x_1,\ldots,x_n,a,b:y_1, \ldots, y_p;m,k,\varepsilon),\frac{r\omega}{64})$.
When $m$ is sufficiently large and $\varepsilon$ is sufficiently
small, we can assume that $A^{\lambda}$ to be $A_k$ and
$B^{\lambda}$ to be $(U^{\lambda})^*B_kU^{\lambda}$ for some unitary
matrix $U$.

For $m$ is sufficiently large and $\varepsilon$ ($\leq
\frac{r\omega}{64}$) is sufficiently small, when
$(T_1,\ldots,T_n,A,B,C)\in \Gamma_R(x_1,\ldots,x_n,\\
a,b,y:y_1,\ldots,y_p;m,k,\varepsilon)$, it follows from Lemma 4 in
\cite{Don} that we can assume that $\|C\|\leq 1$. In addition, it is
clear that $\|AC-CB\|_2\leq \varepsilon$ and
$(T_1,\ldots,T_n,A,B)\in
\Gamma_R(x_1,\ldots,x_n,a,b:y_1,\ldots,y_p;m,k,\varepsilon)$. So
there exist some $\lambda\in\Lambda_k$ and $V\in{\mathcal U}_k$ such
that
$$\|(T_1,\ldots,T_n,A,B)-(VT_1^{\lambda}V^*,\ldots,
VT_n^{\lambda}V^*,VA_kV^*,V(U^{\lambda})^*B_kU^{\lambda}V^*)\|_2\leq
\frac{r\omega}{64}.$$ Hence
$$\|A_kV^*CV-V^*CVU^*B_kU\|_2=\|VA_kV^*C-CV(U^{\lambda})^*B_kU^{\lambda}V^*\|_2\leq
\frac{r\omega}{16}.$$  Then, by Lemma \ref{don, lemma 18}, there
exists a set
$\{Ball(C_{\sigma};\frac{\omega}{4})\}_{\sigma\in\Sigma_k}$ of
$\frac{\omega}{4}$-balls that cover $\{W\in {\mathcal M}_k:
\|W\|\leq 1,\
\|WA_k-(U^{\lambda})^*B_kU^{\lambda}W\|_2<\frac{r\omega}{16}\}$ with
$|\Sigma_k|\leq (\frac{24}{\omega})^{32rk^2}$, i.e., there exists
some $C_{\sigma}$ such that
$$\|V^*CV-C_{\sigma}\|_2=\|C-VC_{\sigma}V^*\|_2\leq \frac{\omega}{4}.$$
Thus
$$\|(T_1,\ldots,T_n,C)-(VT_1^{\lambda}V^*,\ldots,VT_n^{\lambda}V^*,VC_{\sigma}V^*)\|_2\leq
\frac{\omega}{2}.$$ Therefore
$$\nu(\Gamma_R(x_1,\ldots,x_n,y:a,b,y_1,\ldots,y_p);\omega)\leq
|\Lambda_k|\cdot |\Sigma_k|.$$

Hence, we get
\begin{eqnarray*} 0&\leq &
\mathfrak{K}(x_1,\ldots,x_n,y:a,b,y_1,\ldots,y_p;\omega,R)\\
&\leq & \inf_{m\in \Bbb
N,\varepsilon>0}\limsup_{k\rightarrow\infty}\frac{\log(|\Lambda_k||\Sigma_k|)}{-k^2\log\omega}\\
&\leq& \inf_{m\in \Bbb
N,\varepsilon>0}\limsup_{k\rightarrow\infty}\left(\frac{\log(|\Lambda_k|)}{-k^2\log\omega}+\frac{32rk^2(\log
24-\log\omega)}{-k^2\log\omega}\right)\\
&=& \inf_{m\in \Bbb
N,\varepsilon>0}\limsup_{k\rightarrow\infty}\frac{\log(|\Lambda_k|)}{-k^2\log\omega}+32r\frac{\log 24-\log\omega}{-\log\omega}\\
&=&32r\frac{\log 24-\log\omega}{-\log\omega}.
\end{eqnarray*} Since $r$ is an arbitrarily small positive number,
we have
$\mathfrak{K}(x_1,\ldots,x_n,y:a,b,y_1,\ldots,y_p;\omega,R)=0$,
whence, $\mathfrak{K}_2(x_1,\ldots,x_n,y:a,b, y_1,\ldots,y_p)=0$.
Therefore, by theorem \ref{equivalence}, $\mathfrak{K}_3({\mathcal
M})=0$.\hfill $\Box$

\begin{corollary}\label{corollary, x:u=0}Let $\cal M$ be a von Neumann algebra with a faithful trace $\tau$. Suppose $a, x_1, \ldots, x_n$
are elements in $\cal M$ such that $a$ is a normal element without
eigenvalues, and $ax_i=x_ia$ for all $1\leq i\leq n$. Then
$$\mathfrak{K}_3(x_1, \ldots, x_n:u)=0$$
\end{corollary}

Using the similar idea in the proof of Theorem \ref{theorem, without
comon eigenvalues}, we can prove the following theorem.

\begin{theorem}\label{add one unitary}Suppose $\mathcal M$ is a von Neumann algebra with
a faithful trace $\tau$ and ${\mathcal M}=\{{\mathcal N},u\}''$,
where $\mathcal N$ is a von Neumann subalgebra of $\mathcal M$, $u$
is a unitary element of $\mathcal M$. Let $\{v_1,v_2,\ldots\}$ be a
family of Haar unitary elements and $\{w_1,w_2,\ldots\}$ be a family
of unitary elements in $\mathcal N$ such that
$\|v_nu-uw_n\|_2\rightarrow 0$. Then $\mathfrak{K}_3({\mathcal
M})\leq \mathfrak{K}_3({\mathcal N})$.

 In
particular, if $v$ and $w$ are Haar unitary elements in $\mathcal N$
such that $vu=uw$, then $\mathfrak{K}_3({\mathcal M})\leq
\mathfrak{K}_3({\mathcal N})$.
\end{theorem}

%so there is some Borel function $\varphi$ such that $c_0=\varphi(c)$.
%Define another Borel function $\psi$,
%$$\psi(\lambda)=\left\{\begin{array}{ll}-1,&&\lambda\in\sigma_p(b_1)\\
%\varphi(\lambda),&&\lambda\in\sigma(b_1)/\sigma_p(b_1)\\
%0\leg \psi\leq 1&& others\end{array}\right.$$

\section{Applications}

\hspace{1.5em}Suppose $\Lambda$ is an infinite set. An {\it
ultrafilter} $\alpha$ on $\Lambda$ is a collection of subsets of
$\Bbb N$ such that the empty set $\emptyset\not\in\alpha$, $\alpha$
is closed under finite intersections, and, for each subset $A$ of
$\Lambda$, either $A\in\alpha$ or ${\Bbb N}\setminus A\in \alpha$.
One example of an ultrafilter is obtained by choosing an $\iota$ in
$\Lambda$ and letting $\alpha$ be the collection of all subsets of
$\Lambda$ that contain $\iota$. Such an ultrafilter is called {\it
principal}; ultrafilters not of this form are called {\it free}. We
will call an ultrafilter $\alpha$ {\it nontrivial} if it is free and
there exists a decreasing sequence in $\alpha$ whose intersection is
empty. Free ultrafilters on an countable set are always nontrivial.

Suppose $\mathfrak{X}$ is another set, $f:\Lambda\rightarrow
\mathfrak{X}$ is a mapping and $E\subseteq \mathfrak{X}$. we say
that $f(\iota)$ is {\it eventually in $E$} along $\alpha$ if
$f^{-1}(E)=\{\iota\in\Lambda: f(\iota)\in E\}\in \alpha$. If
$\mathfrak{X}$ is a topological space, we say that $f(\iota)$ {\it
converges to $x$} (in $\mathfrak{X}$) along $\alpha$, denoted by
$\lim_{\iota\rightarrow \alpha}f(\iota)=x$, if $f(\iota)$ is
eventually in each neighborhood of $x$. It is well known that if
$\mathfrak{X}$ is a compact Hausdorff space, the
$\lim_{\iota\rightarrow \alpha}f(\iota)$ always exists in
$\mathfrak{X}$ for every $f:\Lambda\rightarrow \mathfrak{X}$ and
every ultrafilter $\alpha$ on $\Lambda$.

Let $\alpha$ be a nontrivial ultrafilter on $\Lambda$. Suppose
${\mathcal M}_\iota$ is a finite von Neumann algebra with a faithful
trace $\tau_\iota$. Let $\prod_\iota{\mathcal M}_\iota$ be the
$l^\infty$-product of the ${\mathcal M}_\iota$'s, ${\mathcal
J}=\{\{x_{\iota}\}: \lim_{\iota\to\alpha}
\tau_{\iota}(x_{\iota}^*x_{\iota})=0\}$. Then define the {\it
ultraproduct} $\prod^\alpha{\mathcal M}_\iota$ of ${\mathcal
M}_\iota$ to be $\prod_{i\in\Bbb I}{\mathcal M}_i/\mathcal J$.

When ${\mathcal M}_\iota=\mathcal M$ for all $\iota$, then
$\prod^\alpha{\mathcal M}_\iota$ is called the {\it ultrapower} of
$\mathcal M$, denoted by ${\mathcal M}^\alpha$.

 Let $\mathcal M$ be a II$_1$ factor with the faithful trace
$\tau$. For every $\varepsilon>0$, and any elements
$x_1,x_2,\ldots,x_n$ in $\cal M$, if there exists a unitary
$u\in\cal M$ with $\tau(u)=0$ such that $\|ux_i-x_iu\|_2\leq
\varepsilon$ for every $i$, then we say that $\mathcal M$ has {\it
property $\Gamma.$}

It is well-known that if $\cal M$ is a II$_1$ factor with property
$\Gamma$, then there exists a central sequence $\{v_n\}$ of Haar
unitary elements in $\mathcal M$, i.e., $\|v_nx-xv_n\|_2\rightarrow
0$ for every $x\in{\mathcal M}$.

If a von Neumann algebra acts on a very large Hilbert space, it may
not contain any nontrivial central sequences, but it may contain a
{\it central net}, i.e., a net $\{x_{\lambda}\}$ in $\cal M$ such
that $\|x_{\lambda}a-ax_{\lambda}\|_2\rightarrow 0$ for every
$a\in\cal M$. Equivalently, $\cal M$ has a central net if and only
if, for every $\varepsilon>0$ and for every finite subset
$F\subseteq {\cal M}$, there is a Haar unitary $u$ such that
$\|ua-au\|_2<\varepsilon$ for every $a\in F$.

\begin{theorem}
\label{ultraproduct}Suppose $\{{\mathcal{M}}_{\lambda }\}_{\lambda
\in
\Lambda }$ is a family of von Neumann algebras such that each $\mathcal{M}%
_{\lambda }$ has a central net of Haar untiaries. Let $\alpha $ be a
nontrivial ultrafilter on $\Lambda $. Then
$$
\mathfrak{K}_{3}(\prod^{\alpha }\mathcal{M}_{\lambda })=0.
$$
\end{theorem}

Proof. Suppose $x_{1},\ldots ,x_{n}$ are any elements in $\prod^{\alpha }%
\mathcal{M}_{\lambda }$, and $x_{i}=\{x_{\lambda }^{i}\}_{\alpha }$. Since $%
\alpha $ is nontrivial, there exists a decreasing sequence $%
\{A_{k}\}_{k=1}^{\infty }$ in $\alpha $ whose intersection is empty, and $%
A_{1}=\Lambda $. If $\lambda \in A_{k}/A_{k-1}$, since ${\mathcal{M}}%
_{\lambda }$ has a central net of Haar unitaries, then there exists
a Haar unitary $u_{\lambda }\in \mathcal{M}_{\lambda }$ such that
$\Vert x_{\lambda
}^{i}u_{\lambda }-u_{\lambda }x_{\lambda }^{i}\Vert _{2}<\frac{1}{k}$ for $%
1\leq i\leq n$. Then $u=\{u_{\lambda }\}_{\alpha }$ defines a Haar
unitary
that commutes with $x_{1},\ldots ,x_{n}$. By Corollary \ref{corollary, x:u=0}%
, $\mathfrak{K}_{2}(x_{1},\ldots ,x_{n}:u)=0.$ Therefore, by Theorem \ref%
{equivalence}, $\mathfrak{K}_{3}\left( \prod^{\alpha
}\mathcal{M}_{\lambda }\right) =0$.\hfill $\Box $

\begin{remark}
Suppose $\cal M$ is a diffuse finite von Neumann algebra with a
faithful trace $\tau$. We can define a numerical invariant $\gamma
\left( \mathcal{M} \right) $ by
$$\gamma \left( \mathcal{M}\right) =\sup_{\begin{array}{c}F\subset
\mbox{ball}\left( \mathcal{M}\right)\\
F\mbox{ finite}\end{array}}\inf_{\begin{array}{c}u\ \mbox{is a}\\
\mbox{Haar unitary}\end{array}}\max_{a\in F}\| au-ua\|.$$
It is clear that $\gamma ( \mathcal{M}) =0$ if and only if $%
\mathcal{M}$ has a central net of Haar unitaries and that $\gamma (
\mathcal{M}) \leq 2$. F. Murray and J. von Neumann \cite{Mv} found a lower bound for $%
\gamma \left( \mathcal{L}_{\mathbb{F}_{n}}\right) $ for $n\geq 2$.
(Of course, they did not use our terminology.) It is not difficult
to modify the proof of Theorem \ref{ultraproduct} to prove
that if $\gamma ({\cal M}_{\lambda})\rightarrow 0$ along the ultrafilter $\alpha $, then $%
\mathfrak{K}_{3}(\prod^{\alpha }\mathcal{M}_{\lambda })=0.$
\end{remark}

\begin{lemma}\cite{Don-Li}
\label{goober}Suppose $\mathcal{M}$ is diffuse, countably generated
and has
a central sequence of Haar unitaries. Then there is a central sequence $%
\{u_{n}\}$ of mutually commuting Haar unitaries in $\mathcal{M}$.
\end{lemma}

\begin{theorem}
If $\mathcal{M}$ is a von Neumann algebra with a central net of Haar
unitaries, then $\mathfrak{K}_{3}({\mathcal{M}})=0$.
\end{theorem}

Proof. Suppose $x_{1},\ldots ,x_{n}\in \mathcal{M}$. Then there is a
sequence $\left\{ u_{n}\right\} $ of Haar unitaries so that if $a\in
\{ x_{1},\ldots ,x_{n},u_{1},u_{2},\ldots\} ,$ then $$\|
au_{n}-u_{n}a\| _{2}\rightarrow 0, $$ i.e., inductively choose
$u_{n}$ so that $\| au_{n}-u_{n}a\|_{2}<1/n$ for $a\in \left\{
x_{1},\ldots ,x_{n},u_{1},u_{2},\ldots u_{n-1}\right\} $. Hence
$\left\{ u_{n}\right\} $
is a central sequence in the von Neumann algebra $\mathcal{N}$ generated by $%
\{ x_{1},\ldots ,x_{n}$, $u_{1},u_{2},\ldots \}$. It follows from
Lemma \ref{goober} that there is a central sequence $\{
v_{n}\}$ of commuting Haar unitaries in $\mathcal{N}$. Then we get $\mathfrak{K}%
_{3}\left( \left\{ v_{1},v_{2},\ldots \right\} ^{\prime \prime
}\right) =0.$ We can choose unitaries $\left\{ w_{1},\ldots
,w_{m}\right\}$ that generate $W^{\ast }\left( x_{1},\ldots
,x_{n}\right),$ and, using Theorem \ref{add one unitary}, we
inductively get  $\mathfrak{K}_{3}\left( \left\{ w_{1},\ldots
,w_{j},v_{1},v_{2},\ldots \right\} ^{\prime \prime }\right) =0$ for
$1\leq j\leq m$. Hence, there exist $y_{1},\ldots ,y_{p}\in
\mathcal{N}$ such that $$ \mathfrak{K}_{2}\left( x_{1},\ldots
,x_{n}:y_{1},\ldots ,y_{p}\right) =0.
$$
Therefore $\mathfrak{K}_{3}({\mathcal{M}})=0.$\hfill $\Box $

\begin{theorem}
Suppose $\mathbb{F}$ is a free group with the standard generating
set $G$ satisfying $|G|\geq 2$ and let $\mathcal{L}_{\mathbb{F}}$ be
the group von Neumann algebra
generated by $\mathbb{F}$. Then $\mathfrak{K}_{3}({\mathcal{L}}_{\mathbb{F}%
})=\infty.$
\end{theorem}

Proof. For any $g\in\Bbb F$, we can view $g$ as a unitary in $\mathcal{L}_{\mathbb{F}}$. Note that $G$ generates ${\mathcal{L}}_{\mathbb{F}}$. Let $%
g_{1},\ldots ,g_{n}\in G$ ($n\geq 2$), and $y_{1},y_{2},\ldots
,y_{N}\in G.$ We will prove that $n\leq \delta \left( g_{1},\ldots
,g_{n}:y_{1},y_{2},\ldots ,y_{N}\right) $. Since $$\delta \left(
g_{1},\ldots ,g_{n}:y_{1},y_{2},\ldots ,y_{N}\right) \leq
1+\mathfrak{K}_{2}\left( g_{1},\ldots ,g_{n}:y_{1},y_{2},\ldots
,y_{N}\right) ,$$ we will conclude that $n-1\leq
\mathfrak{K}_{2}\left( g_{1},\ldots
,g_{n}:y_{1},y_{2},\ldots ,y_{N}\right) .$ From this it follows that  $%
\mathfrak{K}_{3}\left( {\mathcal{L}}_{\mathbb{F}}\right) =\mathfrak{K}%
_{3}\left( G\right) =\infty $.

It follows from Theorem 13 in \cite{Don} that when we compute
$\delta $ we
can replace the $\Gamma _{R}$-sets with the set of unitary elements in the $%
\Gamma _{R}$-sets. Let $$\Omega _{m,k,\varepsilon
}=\mathcal{U}_{k}^{n+N}\cap \Gamma_R (g_{1},\ldots
,g_{n},y_{1},y_{2},\ldots ,y_{N};m,k,\varepsilon ),$$  $$\Delta
_{m,k,\varepsilon }=\mathcal{U}_{k}^{n}\cap \Gamma _{R}(g_{1},\ldots
,g_{n}:y_{1},y_{2},\ldots ,y_{N};m,k,\varepsilon ).$$ If
$(U_{1},\ldots
,U_{n},W_{1},\ldots ,W_{N})\in \Omega _{m,k,\varepsilon }$, then  $%
(U_{1},\ldots ,U_{n})\in \Delta _{m,k,\varepsilon }$.  So $$\Omega _{m,k,\varepsilon }\subseteq \Delta _{m,k,\varepsilon }\times {\mathcal{%
U}}_{k}^{N}.
$$
Let $\mu _{k}$ denote Haar measure on ${\mathcal{U}}_{k}$, $\mu
_{k}^{n}$ denote the corresponding product measure on
${\mathcal{U}}_{k}^{n}$. Then $$ \mu _{k}^{n+N}\left( \Omega
_{m,k,\varepsilon }\right) \leq \mu
_{k}^{n}(\Delta _{m,k,\varepsilon })\cdot \mu _{k}^{N}\left( {\mathcal{U}}%
_{k}^{N}\right) \leq \mu _{k}^{n}(\Delta _{m,k,\varepsilon })\leq 1.
$$
We know from Theorem 3.9 in \cite{V02} that $\mu _{k}^{n+N}\left(
\Omega _{m,k,\varepsilon }\right) \rightarrow 1$ and thus $\mu
_{k}^{n}(\Delta _{m,k,\varepsilon })\rightarrow 1 $ as $k\rightarrow
\infty $. It follows that $$ \mu _{k}^{n}(\Delta _{m,k,\varepsilon
})\leq \nu _{2}(\Delta _{m,k,\varepsilon },\omega )\mu
_{k}^{n}(ball(1,\omega ))\leq \nu _{2}(\Delta _{m,k,\varepsilon
},\omega )(\omega )^{nk^{2}},
$$
so $\nu _{2}(\Delta _{m,k,\varepsilon },\omega )\geq \mu
_{k}^{n}(\Delta _{m,k,\varepsilon })(\frac{1}{\omega })^{nk^{2}}$.

By Lemma 1 in \cite{Don-Shen}, $$ \mathfrak{K}_{2}(g_1, \ldots,
g_n:y_1, \ldots, y_N)\geq  \delta (g_1, \ldots, g_n:y_1, \ldots,
y_N) -1.$$ Note that \begin{eqnarray*} \delta (g_1, \ldots, g_n:y_1,
\ldots, y_N)&=&-1+\limsup_{\omega \rightarrow
0^{+}}\inf_{m,\varepsilon }\limsup_{k\rightarrow \infty }\frac{\log
\left( \nu _{2}(\Delta
_{m,k,\varepsilon },\omega )\right) }{-k^{2}\log \omega }\\
& \geq &-1+\limsup_{\omega \rightarrow 0^{+}}\inf_{m,\varepsilon
}\limsup_{k\rightarrow \infty }\frac{\log \left( \mu _{k}^{n}(\Delta
_{m,k,\varepsilon })\right) -nk^{2}\log \omega }{-k^{2}\log \omega
}.\end{eqnarray*} Since $\mu _{k}^{n}(\Delta _{m,k,\varepsilon
})\rightarrow 1$, we have $$\mathfrak{K}_{2}(g_1, \ldots, g_n:y_1,
\ldots, y_N)\geq n-1.$$ Thus
$\mathfrak{K}_{3}({\mathcal{L}}_{G})=\infty .$\hfill $\Box $

\begin{theorem}Suppose $\mathcal M$ is a von Neumann algebra with
a faithful trace $\tau$, ${\mathcal N}_1$ and ${\mathcal N}_2$ are
mutually commuting diffuse subalgebras of $\mathcal M$. Then
$\mathfrak{K}_3(W^*({\cal N}_1\cup{\cal N}_2))=0$.
\end{theorem}

Proof. Since ${\mathcal N}_1$ and ${\mathcal N}_2$ are diffuse, we
can assume that ${\mathcal N}_1=\{u_{\lambda}:\lambda\in\Lambda\}''$
and ${\mathcal N}_2=\{v_{\sigma}:\sigma\in\Sigma\}''$ where
$u_{\lambda}, v_{\sigma}$ are all Haar unitaries. For any finite
subset $E$ of $\Lambda$ and finite subset $F$ of $\Sigma$, let
${\mathcal M_{E,F}}=\{u_{\lambda}, v_{\sigma}: \lambda\in E,
\sigma\in F\}''$. Then $\mathfrak{K}_3\left({\mathcal
M_{E,F}}\right)=0$ by Theorem \ref{add one unitary}. Let
$${\mathcal S}=\bigcup\left\{{\mathcal M_{E,F}}: E\ \mbox{is a
finite subset of}\ \Lambda, F\ \mbox{is a finite subset of}\
\Sigma\right\}.$$ It is clear that ${\mathcal S}$ generates
$W^*({\cal N}_1\cup{\cal N}_2)$. By Corollary \ref{corollary,
increasingly directed}, $\mathfrak{K}_3({\mathcal S})=0$. Thus
$\mathfrak{K}_3(W^*({\cal N}_1\cup{\cal N}_2))=0$ by Corollary
\ref{corollary, k3 is independent on generatoer}.\hfill $\Box$

\begin{corollary}If ${\mathcal N}_1$ and
${\mathcal N}_2$ are diffuse von Neumann algebras, then
$\mathfrak{K}_3({\mathcal N}_1\otimes {\mathcal N}_2)=0$.
\end{corollary}

The following corollary was proved by S. Popa \cite{popa2} and L. Ge
\cite{Ge2}.
\begin{corollary}\label{corollary, free group factor k_3=infty}If ${\mathcal M}$ is a II$_1$ factor and $\mathfrak{K}_3({\cal
M})=\infty$, then $\cal M$ is prime. In particular,
${\mathcal{L}}_{\mathbb{F}}$ is prime for every free group $\Bbb F$
with the standard generating set $G$ satisfying $|G|\geq 2$..
\end{corollary}

Let $\mathcal M$ be a von Neumann algebra with a faithful trace,
${\cal U(M)}$ be the set of all unitary elements in $\cal M$. For
any subset ${\cal S}\subseteq {\cal M}$, define $${\cal N}({\cal
S})=\{u\in {\cal U(M)}: u{\cal S}u^*\subseteq {\cal S}\}'',$$ and
\begin{eqnarray*}{\cal I}({\cal S})&=&W^*(\{y\in{\cal M}: \exists\  \mbox{two normal elements}\ a,b \\
&&\mbox{without common eigenvalues such that}\ ay=yb\neq
0\})\end{eqnarray*} Suppose $\mathcal A$ is a diffuse von Neumann
subalgebra of $\mathcal M$, and $\alpha$ is an ordinal. Then define
$${\mathcal N}_{\alpha}({\mathcal
A})=\left\{\begin{array}{lll}\mathcal A& \alpha=0\\
\left(\bigcup_{\beta<\alpha}{\mathcal N}_{\beta}(\mathcal
A)\right)''& \alpha\ \mbox{is a limit ordinal}\\
{\cal N}\left({\cal N}_{\beta}({\cal A})\right) & \mbox{if}\
\alpha=\beta+1,\end{array}\right.$$ and
$${\mathcal I}_{\alpha}({\mathcal
A})=\left\{\begin{array}{lll}\mathcal A& \alpha=0\\
\left(\bigcup_{\beta<\alpha}{\mathcal I}_{\beta}(\mathcal
A)\right)''& \alpha\ \mbox{is a limit ordinal}\\
{\cal I}\left({\cal I}_{\beta}({\cal A})\right) & \mbox{if}\
\alpha=\beta+1.\end{array}\right.$$

The following theorem is a easy consequence of Theorem \ref{add one
unitary}, Theorem \ref{theorem, without comon eigenvalues} and
Corollary \ref{corollary, increasingly directed}.

\begin{theorem}\label{thereom, normalizer}Let $\mathcal M$ be a von Neumann algebra with
faithful trace $\tau$ and $\mathcal A$ be a diffuse subalgebra of
$\mathcal M$ with $\mathfrak{K}_3({\mathcal A})=0.$ Then for any
ordinal $\alpha$, $$\mathfrak{K}_3\left({\mathcal
N}_{\alpha}({\mathcal A})\right)=\mathfrak{K}_3\left({\mathcal
I}_{\alpha}({\mathcal A})\right)=0.$$
\end{theorem}

\section{Applications to group theory}

\hspace{1.5em}Suppose $\cal M$ is a von Neumann algebra with a
faithful trace $\tau$ and $u$ is a Haar unitary in $\cal M$. Define
$${\cal N}_u=\left(\cup\{ {\cal N}: {\cal N}\subseteq {\cal M}, u\in
{\cal N}, \mathfrak{K}_3({\cal N})=0\}\right)''$$ to be the von
Neumann subalgebra generated by the union of those subalgebras $\cal
N$ containing $u$ such that $\mathfrak{K}_3({\cal N})=0$. It follows
from Theorem \ref{theorem, intersection contains Haar unitary} that
$\mathfrak{K}_3\left({\cal N}_u\right)=0.$ Therefore ${\cal N}_u$ is
the unique largest subalgebra containing $u$ with
$\mathfrak{K}_3({\cal N}_u)=0$.

Suppose $\cal B$ is a diffuse von Neumann subalgebra of $\cal M$
with $\mathfrak{K}_3({\cal B})=0$. Since $\cal B$ is diffuse, there
exists a Haar unitary $u\in \cal B$. It is clear that ${\cal N}_u$
is the largest subalgebra of $\cal M$ that contains $\cal B$ whose
$\mathfrak{K}_3$ is 0.

If $u, v$ are two Haar unitaries in $\cal M$ such that ${\cal
N}_u\cap {\cal N}_v$ is diffuse, then, by Theorem \ref{theorem,
intersection contains Haar unitary},
$\mathfrak{K}_3\left(\left({\cal N}_u\cup {\cal
N}_v\right)''\right)=0.$ Therefore ${\cal N}_u={\cal N}_v$ by the
maximality of ${\cal N}_u$ and ${\cal N}_v$. Therefore, if ${\cal
N}_u\neq {\cal N}_v$, then ${\cal N}_u\cap {\cal N}_v$ is not
diffuse.

By Theorem \ref{thereom, normalizer}, it is clear that, for any Haar
unitary $u\in \cal M$, the normalizer
 of ${\cal N}_u$ is ${\cal N}_u$.

The above ideas have an interesting interpretation in discrete
groups. Let $G$ be a group and ${\cal L}_G$ be the corresponding
group von Neumann algebra.  For any $g\in G$, we can view $g$ as a
unitary in ${\cal L}_G$. If $g$ is an element in $G$ with infinite
order, then $g$ is a Haar unitary in ${\cal L}_G$. If $H$ is a
subgroup of $G\subseteq {\cal L}_G$, then $H''\cong {\cal L}_H.$

Suppose $G$ is a discrete torsion-free group (i.e., the only element
of finite order is identity) and ${\mathcal L}_{G}$ can be embedded
into an ultrapower of the hyperfinite II$_1$ factor. F. Radulescu
(see Proposition 2.5 in \cite{R}) proved that this is equivalent to
$G$ being algebraically embeddable in the unitary group of such an
ultrapower. For any $g\in G\setminus \{e\}$, let $H_g$ be the
subgroup generated by the set $\{H\leq G: g\in H,
\mathfrak{K}_3({\cal L}_{H})=0\}.$

By Theorem \ref{theorem, intersection contains Haar unitary} and
Theorem \ref{add one unitary}, we get the following theorem.
\begin{theorem}Let $G$ be a torsion-free group with the unit $e$.
The following statements hold:

(1) $\{H_g\setminus\{e\}: g\in G\setminus\{e\}\}$ is a partition of
$G\setminus\{e\}$,

(2) for any $g\in G\setminus \{e\}$, if $hH_gh^{-1}\subseteq H_g$,
then $h\in H_g$,

(3) for any $g\in G\setminus\{e\}$, if $h\in G\setminus H_g$, then
$hH_g\cap H_gh=\{h\}$.

\end{theorem}

We call $\{H_g: g\in G/\{e\}\}$ the {\it
$\mathfrak{K}_3$-decomposition} of the torsion-free group $G$.
Determining the $\mathfrak{K}_3$-decomposition of a particular group
involves both algebra and the theory of von Neumann algebras. In
\cite {Don-Shen}, there are many examples of groups $G$ whose
$\mathfrak{K}_3$-decomposition is $\{G\}$. We now provide a
particular example.

\begin{proposition}If $G$ is a free group, then, for any $g\in G\setminus \{e\}$,
$H_g$ is a maximal cyclic subgroup of $G$ containing $g$. In
particular, if $g$ is one of the free generators, then $H_g$ is the
subgroup generated by $g$.
\end{proposition}

Proof. Since every subgroup of free group is free, $H_g$ is free. By
definition of $H_g$ we see that $\mathfrak{K}_3({\cal L}_{H_g})=0$.
Therefore, by
 Corollary \ref{corollary, free group factor k_3=infty}, $H_g$ cannot have more than one generator. That implies $H_g$ is cyclic. By Theorem
 \ref{add one unitary}, $H_g$ must be a maximal abelian subgroup in
 $G$.\hfill $\Box$
\vspace{0.2cm}

Let ${\Bbb F}_2$ be a free group generated by two standard
generators $u,v$. We know that, in ${\Bbb F}_2$, $H_u$ is the
subgroup generated by $u$. This naturally raises the question.

\vspace{0.3cm}

\noindent \textbf{Question 1}.  In ${\cal L}_{{\Bbb F}_2}$, is
${\cal N}_u=W^*(u)$? In other words, is $W^*(u)$ a maximal
subalgebra of ${\mathcal L}_{{\Bbb F}_2}$ whose $\mathfrak{K}_3$ is
$0$?

 S. Popa \cite{popa3} proved that $W^*(u)$ is maximal injective. An
 affirmative answer to the question above would imply Popa's result,
 since $\mathfrak{K}_3({\cal M})=0$ whenever ${\cal M}$ is
 injective. This means that answering the question above is likely
 to be difficult. However, there are natural subquestions based on
 Theorem \ref{theorem, intersection contains Haar unitary} and Theorem
 \ref{theorem, without comon eigenvalues}, respectively.

\vspace{0.3cm}

\noindent \textbf{Question 1a} If $\cal M$ is a subalgebra of
${\mathcal L}_{{\Bbb F}_2}$ with $\mathfrak{K}_3({\cal M})=0$, and
${\cal M}\cap W^*(u)$ is diffuse, then must we have ${\cal
M}\subseteq W^*(u)$?

\vspace{0.3cm}

\noindent \textbf{Question 1b } If $y\in {\mathcal L}_{{\Bbb F}_2}$,
and $a, b\in W^*(u)$ without common eigenvalues, such that,
$ya=by\neq 0$, then must y be in $W^*(u)$?

\vspace{0.2cm}

We can give a partial solution to Question 1b by showing that if $w$
is a unitary in ${\mathcal L}_{{\Bbb F}_2}$ that conjugates a Haar
unitary in $W^*(u)$ into $W^*(u)$, then $w\in W^*(u)$.

 Suppose that $\cal M$ and $\cal N$ are von Neumann algebras and
${\cal N}\subseteq {\cal M}$. By a {\it conditional expectation}
from $\cal M$ onto $\cal N$, we mean a positive linear mapping $E:
{\cal M}\rightarrow {\cal N}$ such that

(1) $E(I)=I$,

(2) $E(x_1yx_2)=x_1E(y)x_2$, for any $x_1, x_2\in\cal N$ and
$y\in\cal M$.

 Define ${\mathcal
A}\bot{\mathcal B}$ to be $\tau(ab)=\tau(a)\tau(b)$ for any
$a\in\cal A$ and $b\in\cal B$.

\begin{theorem} Let $u,v$ be standard generators of ${\cal L}_{{\Bbb F}_2}$. If $\mathcal B$ is a diffuse von Neumann subalgebra of
$W^*(u)$, then $$\{w:\ w\ \mbox{is a unitary in}\ {\cal L}_{{\Bbb
F}_2}, w^*{\cal B}w\subseteq W^*(u)\}\subseteq W^*(u).$$
\end{theorem}

Proof. Suppose $\mathcal B$ is a diffuse von Neumann subalgebra of
$W^*(u)$. Define $${\cal N}=\{w:\ w\ \mbox{is a unitary in}\ {\cal
L}_{{\Bbb F}_2}, w^*{\cal B}w\subseteq {W^*(u)}\}.$$ It is
sufficient to prove $W^*(u)^{\perp}\subseteq {\cal N}^{\perp}$. The
proof is a modification of Lemma 2.5 in \cite{popa2}.

Let $g$ be en element in $W^*(u)^{\perp}$. It follows that
$gW^*(u)g^*\bot W^*(u)$.

Since $\cal B$ has no atoms, then, for any given $\varepsilon>0$,
there exists an orthogonal family of projections $e_1,\ldots,e_n$ in
$\cal B$ such that $\tau(e_i)<\varepsilon$ for $1\leq i\leq n$. Let
${\mathcal A}_{\varepsilon}$ be the von Neumann subalgebra generated
by $e_1,\ldots,e_n$, $\tau$ be the unique trace on ${\mathcal
L}_{{\Bbb F}_2}$, and $E_{{\cal A}_{\varepsilon}'\cap {\cal
L}_{{\Bbb F}_2}}$ be the unique $\tau$-preserving conditional
expectation from ${\mathcal L}_{{\Bbb F}_2}$ onto ${\cal
A}_{\varepsilon}'\cap {\cal L}_{{\Bbb F}_2}$ (i.e., $\tau\circ
E_{{\cal A}_{\varepsilon}'\cap {\cal L}_{{\Bbb F}_2}}=\tau$).

For any $w\in {\cal N}$, $g{\cal A}_{\varepsilon}g^*\bot w{\cal
A}_{\varepsilon}w^*$ since $g{\cal A}_{\varepsilon}g^*\bot W^*(u)$
and $w{\cal A}_{\varepsilon}w^*\subseteq w{\cal B}w^*\subseteq
W^*(u)$. Therefore, for $1\leq i\leq n$,
$$\tau(w^*ge_ig^*we_i)=\tau(ge_ig^*we_iw^*)=\tau(ge_ig^*)\tau(we_iw^*)=\tau(e_i)^2.$$
Summing up over $i$, we get
\begin{eqnarray*}|\tau(wg)|^2&\leq & \|E_{{\cal A}_{\varepsilon}'\cap {\cal
L}_{{\Bbb F}_2}}(wg)\|_2^2\\
&=&\|\sum_{i}e_iwge_i\|_2^2=\sum_{i}\|e_iwge_i\|_2^2\\
&=&\sum_{i}\tau(wge_ig^*w^*e_i)=\sum_i\tau(e_i)^2\\
&\leq& (\max_j\tau(e_j))\sum_i\tau(e_i)<\varepsilon.
\end{eqnarray*}
Since $\varepsilon>0$ is arbitrarily small, $\tau(wg)=0$. Therefore
for any $x\in {\cal N}''$, $\tau(x)=0$. Thus $g\bot {\cal N}''$.
\hfill$\Box$

\end{document}